\documentstyle[12pt,amsfonts,leqno,amscd]{amsart}
\newtheorem{Theorem}{Theorem}[section]
\newtheorem{Definition}[Theorem]{Definition}
\newtheorem{Proposition}[Theorem]{Proposition}

\newtheorem{Lemma}[Theorem]{Lemma}
\newtheorem{Corollary}[Theorem]{Corollary}
\theoremstyle{remark}
\newtheorem{Example}[Theorem]{Example}

\def\il{\int\limits_}

\def\ovr{\overline}

\def\om{\omega}

\def\al{\alpha}

\def\th{\theta}

\def\Gm{\Gamma}

\def\sm{\setminus}
\def\sbs{\subset}

\def\wtl{\widetilde}

\def\be{\begin{enumerate}}
\def\ee{\end{enumerate}}
\def\bT{\begin{Theorem}}
\def\eT{\end{Theorem}}
\def\bP{\begin{Proposition}}
\def\eP{\end{Proposition}}
\def\bD{\begin{Definition}}
\def\eD{\end{Definition}}
\def\bE{\begin{Example}}
\def\eE{\end{Example}}
\def\bL{\begin{Lemma}}
\def\eL{\end{Lemma}}
\def\bC{\begin{Corollary}}
\def\eC{\end{Corollary}}

\newcommand{\ra}{\rightarrow}
\begin{document}
\title{Quasianalyticity and pluripolarity}
\author{ Dan Coman, Norman Levenberg and Evgeny A. Poletsky}\date{November 2002}
\keywords{Quasianalytic functions, pluripolar sets, pluripotential
theory} \subjclass[2000]{ Primary: 26E10, 32U20; secondary: 32U35,
32U15, 32U05}
\thanks{D. Coman and E.A. Poletsky were supported by NSF grants}
\address{D. Coman and E. A. Poletsky: dcoman@@syr.edu, eapolets@@syr.edu,
Department of Mathematics,  215 Carnegie Hall, Syracuse
University, Syracuse, NY 13244, USA}
\address{N. Levenberg: levenber@@math.auckland.ac.nz, Department of
Mathematics, University of Auckland, Private Bag 92019 Auckland,
New Zealand}
\begin{abstract}
We show that the graph
$$\Gm_f=\{(z,f(z))\in{\Bbb C}^2:\,z\in S\}$$
in ${\Bbb C}^2$ of a function $f$ on the unit circle $S$ which is
either continuous and quasianalytic in the sense of Bernstein or
$C^\infty$ and quasianalytic in the sense of Denjoy is pluripolar.
\end{abstract}
\maketitle
\section{Introduction}
A set $K$ in ${\Bbb C}^n$ is called {\it pluripolar} if there is a
plurisubharmonic (psh) function $u\not \equiv -\infty$ which is
equal to $-\infty$ on $K$. As an example, if $K$ is an analytic
set given by an equation $h=0$, where $h$ is a holomorphic
function, then $u=\log|h|$ is psh and equal to $-\infty$ on $K$.
\par Let $f$ be a function on the unit circle $S$ and let
$$\Gm_f=\{(z,f(z))\in{\Bbb C}^2:\,z\in S\}$$ be the graph of $f$
in ${\Bbb C}^2$.
\par The set $\Gm_f$ is always pluripolar when $f$ is a {\it
real-analytic}
function. In \cite{DF}, Diederich and Forn\ae ss give an example
of a $C^{\infty}$ function $f$ with non-pluripolar graph in ${\Bbb
C}^2$. The paper \cite{LMP} contains an example of a holomorphic
function $f$ on the unit disk $U$, continuous up to the boundary,
such that the graph of $f$ over $S$  is not pluripolar as a subset
of ${\Bbb C}^2$. Thus {\it a priori } the pluripolarity of graphs of
functions on $S$ is indeterminate.
\par In this paper we prove that graphs of {\it quasianalytic}
functions  are still pluripolar (all necessary definitions can be
found in the next section). More precisely,
\begin{Theorem}\label{T:bq}
If $f:S\ra{\Bbb C}$ is quasianalytic in the sense of Bernstein or
Denjoy, then the set $\Gamma_f\subset{\Bbb C}^2$ is pluripolar.
\end{Theorem}
We also prove
\begin{Theorem}\label{T:dq}
Let $f:S\ra{\Bbb C}$ be a $C^\infty$ function. If $f$ belongs to a
Gevrey class $G\{L_j\}$, where $L_j\geq j$ satisfies $L_j=o(j^2)$
as $j\ra\infty$, then the set $\Gamma_f$ is pluripolar.
\end{Theorem}
The function $f$ in the example of \cite{DF} is in a Gevrey class
$G\{L_j\}$ with $L_j=C^j$; thus, not all functions in a Gevrey
class have pluripolar graphs. On the other hand, taking
$L_j=j^{3/2}$, we see that there are functions in a Gevrey class
that always have pluripolar graphs but which are {\it not}
quasianalytic; i.e., in the $C^{\infty}$ category,
``quasianalytic'' implies ``pluripolar graph'' but not conversely.
\par Since any continuous function is the difference of two
continuous functions which are quasianalytic in the sense of
Bernstein, we arrive at a surprising result.
\bC \label{C:bq} Any continuous function defined on $S$ is the difference
of two continuous functions with pluripolar graphs.\eC
\par According to \cite{M}, any $C^\infty$ function on the unit
circle is the difference of two $C^\infty$ functions which are
quasianalytic in the sense of Denjoy. Hence we have the $C^\infty$
analogue of Corollary \ref{C:bq}.
\bC Any $C^\infty$ function on the unit
circle is the difference of two $C^\infty$ functions with
pluripolar graphs.
\eC
\par In the next section we recall all necessary facts and
definitions. In Section \ref{S:bq} we prove that quasianalytic
functions in the sense of Bernstein have {\it negligible} and,
consequently, pluripolar graphs. To deal with $C^\infty$ functions
we establish in Section \ref{S:aux} a criterion (Theorem
\ref{C:cpp}) for pluripolarity: if the total Monge--Amp\`ere
masses of a sequence of multipole pluricomplex Green functions are
uniformly bounded from above, then the set where the limit of
this sequence is equal to $-\infty$ is pluripolar. This criterion
implies Corollary \ref{P:annuli} which allows us to verify the
pluripolarity of a set by constructing sequences of holomorphic
mappings on annuli with appropriate radii. Theorem \ref{T:bq} for
quasianalytic functions in the sense of Denjoy and Theorem
\ref{T:dq} are proved in Section \ref{S:dq}. For this we
interpolate $f$ at the $n$-th roots of unity and an arbitrary
point $z_0\in S$ by Lagrange trigonometric polynomials $L_n$ and
show that these interpolants are uniformly bounded on annuli
$A(t_n)=\{1/t_n<|z|<t_n\}$, where
$$\limsup_{n\ra\infty}\sqrt{n}\,\log  t_n=\infty.$$
This formula is exactly what we need to apply
Corollary  \ref{P:annuli}.
\par We would like to thank Al Taylor for introducing us to the notion of
quasianalytic functions.
\section{Basic definitions and facts}\label{S:bdf}
\par  For a continuous function $f$ on $S$, we consider the approximation
numbers $E_n(f)=\inf \|f-p_n\|_S$, where $p_n$ runs over all
trigonometric polynomials of degree at most $n$, i.e.,
$$p_n(z)=\sum_{k=-n}^nc_kz^k,\;z\in S,$$
and $\|\cdot\|_S$ is the uniform norm. A continuous function
$f:S\to {\Bbb C}$ is called {\it quasianalytic in the sense of
Bernstein} if
$$\liminf_{n\to \infty} E_n^{1/n}(f)<1.$$
This class contains some continuous functions which are nowhere
differentiable. We refer to \cite{Ti} for further details.

\par Let $f:S\ra{\Bbb C}$ be a $C^\infty$ function with Fourier
series given by
\begin{equation}\label{e:Fourier}
f(z)=\sum_{k=-\infty}^\infty c_kz^k,\;z=e^{i\theta}.
\end{equation}
The $L^2$ norm $M_j(f)$ of the $j$-th derivative $\widetilde
f^{(j)}$, where $\widetilde f(\theta)=f(e^{i\theta})$, is given by
\begin{equation}\label{e:deriv}
M_j^2(f)=\frac{1}{2\pi}\int_0^{2\pi}|\widetilde
f^{(j)}(\theta)|^2d\theta= \sum_{k=-\infty}^\infty k^{2j}|c_k|^2.
\end{equation}
The sequence $\{M_j(f)\}$ is increasing and logarithmically
convex. Following \cite{Ka}, we consider the associated function
$\tau_f(r)$ (associated to the sequence $\{M_j(f)\}$) defined by
$$\tau_f(r)=\inf_{j\geq0}\frac{M_j(f)}{r^j}\;,\;r>0.$$
We can assume that $M_j(f)$ increases faster than $R^j$ for any
$R>0$, or else $f$ is a trigonometric polynomial. Then
$\tau_f(r)>0$ is a decreasing function with
$\lim_{r\ra\infty}\tau_f(r)=0$. The function $-\log\tau_f(r)$
is a convex, increasing function of $\log r$.
\par Now we recall some facts about quasianalytic and Gevrey classes
of smooth functions (see \cite{Ka} and \cite{KP}). Given an
increasing sequence $\{M_j\}$ which is logarithmically convex, the
class $C^\#\{M_j\}$ consists of all smooth functions $f:S\ra{\Bbb
C}$ satisfying the estimate $M_j(f)\leq R^jM_j$ for all $j$ with a
constant $R$ depending on $f$. The class $C^\#\{M_j\}$ is called
quasianalytic if every function in $C^\#\{M_j\}$ which vanishes to
infinite order at some point in $S$ must be identically equal to
0. Let $\tau(r)=\inf_{j\geq0}(M_j/r^j)$ be the associated function
to the sequence $M_j$. The Denjoy-Carleman theorem states that the
class $C^\#\{M_j\}$ is quasianalytic if and only if
\begin{equation}\label{Q:def}
\int_1^\infty\frac{\log\tau(r)}{1+r^2}\;dr=-\infty.
\end{equation}
A smooth function
$f$ is called {\it quasianalytic in the sense of Denjoy} if the class
$C^\#\{M_j(f)\}$ is quasianalytic.

\par Let $\{L_j\}$ be a sequence such that $j\leq L_j\leq CL_{j-1}$
holds for all $j>0$ with some constant $C>0$ independent of $j$.
The {\it Gevrey class} $G\{L_j\}$ consists of all smooth functions $f$
which satisfy $M_j(f)\leq(C'L_j)^j$ for all $j$, where the constant
$C'>0$ depends on $f$. The Gevrey class $G\{L_j\}$ is
quasianalytic if and only if $\sum_{j=0}^\infty(1/L_j)=\infty$.
\section{Quasianalytic functions in the sense of
Bernstein}\label{S:bq}
\par We will consider a trigonometric polynomial
$p(e^{i\theta})=\sum_{k=-n}^nc_ke^{ik\theta}$ as the restriction
to $S$ of the rational function $p(z)=\sum_{k=-n}^nc_kz^k$. The
following lemma is a simple version of a Bernstein--Walsh
inequality for the punctured plane ${\Bbb C}^*={\Bbb C}\sm\{0\}$.
\bL\label{L:betp} Let $p$ be a trigonometric polynomial on $S$ of
degree $n$. Then $|p(z)|\le\|p\|_Se^{nV(z)}$ for all $z\in{\Bbb
C}^*$ where $V(z)=\left|\log|z|\right|$.\eL
\begin{pf} Assume that $\|p\|_S=1$. The function $u(z)=(\log|p(z)|)/n$
is subharmonic when $z\ne0$. Clearly $v(z)=u(z)-V(z)$ is
subharmonic on $U\sm\{0\}$ and on ${\Bbb C}\sm\ovr U$. Also $v$ is
bounded near zero and near infinity. Since $v\le0$ on $S$, we see
that $v\le0$ when $z\ne0$. Thus $|p(z)|\le e^{nV(z)}$ and the
lemma is proved.
\end{pf}
\par Now we can prove our first result regarding the pluripolarity of
graphs.
\bT\label{T:qab} If a function $f:\,S\to{\Bbb C}$ is
quasianalytic in the sense of Bernstein, then $\Gm_f$ is
pluripolar in ${\Bbb C}^2$.\eT
\begin{pf}
We can find $c<1$, a sequence of positive integers $\{n_k\}$, and
a corresponding sequence of trigonometric polynomials
$\{p_{n_k}\}$ with
$$E_{n_k}(f)=\|f-p_{n_k}\|_S \leq c^{n_k}.$$
Without loss of generality, we may assume $\|f\|_S\leq 1/2$ so
that $\|p_{n_k}\|_S\leq 1$; by Lemma \ref{L:betp}
$|p_{n_k}(z)|\leq e^{n_kV(z)}$ for each $k$ and for all $z\in{\Bbb
C}^*$.
\par Define functions
$v_k(z,w):={1\over n_k} \log {|w-p_{n_k}(z)|}$ on ${\Bbb
C}^*\times{\Bbb C}$. Then
$$v_k(z,w)\leq \max\left\{V(z),{1\over n_k} \log|w|\right\}+
\frac{\log 2}{n_k}.$$ Also, if $(z,w)\in\Gm_f$, i.e., $|z|=1$ and
$w=f(z)$, then $v_k(z,w)\le\log c$. Let $v(z,w)=\sup_k v_k(z,w)$.
The function $v$ is bounded above on compacta in ${\Bbb
C}^*\times{\Bbb C}$. Moreover, $v\le\log c$ on $\Gm_f$ and since
$p_{n_k}(z)\to f(z)$ as $k\to\infty$, $v\ge0$ on $(S\times{\Bbb
C})\sm\Gm_f$. Thus the function $v$ is not upper semicontinuous on
$\Gm_f$, i.e., $$\Gm_f\subset \{(z,w)\in {\Bbb C}^*\times{\Bbb
C}:v(z,w) <v^*(z,w)=\limsup_{(z',w')\to (z,w)}v(z',w')\}$$ is a
negligible set. By \cite{BT}, negligible sets are pluripolar and
our theorem is proved.
\end{pf}
\par
\section{A criterion for pluripolarity}\label{S:aux}
\par To avoid notational confusion, we temporarily utilize ${\bf z}$
for a point ${\bf z}=(z_1,...,z_n)\in {\Bbb C}^n$ and write
$\|{\bf z}\|^2=|z_1|^2+\cdots +|z_n|^2$. Let $D$ be a strongly
pseudoconvex domain in ${\Bbb C}^n$ with a strongly psh defining
function $\rho\in C^2(\ovr D)$. In \cite{Sa} Sadullaev proved that
for any set $E\sbs D_r=\{{\bf z}\in D:\rho ({\bf z})\le r\}$,
$r<0$, there are positive constants $\al(r)$ and $\beta(r)$
depending only on $r$ such that
\begin{equation}\label{e:si}
\al(r)C^*(E)\le\il D|\om^*({\bf
z},E,D)|(dd^c\rho)^n\le\beta(r)\left[C^*(E)\right]^{1/n}.
\end{equation}
Here $\om^*({\bf z},E,D)=\limsup_{\zeta \to {\bf z}}\om
(\zeta,E,D)$ is the upper semicontinuous regularization of the
relative extremal function
$$\om (\zeta,E,D)=\sup \{u(\zeta): u \ \hbox{psh in} \ D, \ u\leq 0, \
u|_E\leq -1\}$$ of the set $E$ in $D$, and $C^*(E)$ is the outer
capacity of $E$  relative to $D$:
$$C^*(E)=\inf_V\sup_u\il V(dd^cu)^n,$$
where $u$ runs over all psh functions $u$ on $D$ such that $-1\le
u\le 0$ and $V$ is any open set in $D$ containing $E$ (for a
discussion of the Monge-Amp\`ere operator, $(dd^c(\cdot))^n$, we
refer the reader to \cite{BT}).
\par We also
consider weighted multipole pluricomplex Green functions $G({\bf
z})=G({\bf z};a;\alpha)$ where $a=\{a_1,...,a_m\}$ are points in
$D$ and $\alpha =\{\alpha_1,...,\alpha_m\}$ are positive numbers:
$$G({\bf z})=\sup \{u({\bf z}): u \ \hbox{psh in} \ D, \ u\leq 0, $$
$$\ u({\bf z})-\alpha_j\log {||{\bf z}-a_j||}=O(1), \ {\bf z}\to a_j, \
j=1,...,m\}.$$ It is known that $G$ is continuous and psh in $D$,
$(dd^cG)^N =0$ on $D\setminus \{a_1,...,a_m\}$, $G=0$ on $\partial
D$, and $G({\bf z})-\alpha_j\log {||{\bf z}-a_j||}=O(1), \ {\bf
z}\to a_j, \ j=1,...,m$. We refer to \cite{De} and \cite{Lel} for
further properties of these functions. We have the following
lemma.
\bL\label{L:ma} Suppose that $D=\{\rho<0\}$ is a strongly
pseudoconvex domain in ${\Bbb C}^n$ as above. If $G({\bf z})$ is a
multipole Green function with poles in the set $D_r$, $r<0$, then
there is a positive constant $b(r)$ depending only on $r$ such
that
$$\il D|G|(dd^c\rho)^n\le b(r)\left(\il
D(dd^cG)^n\right)^{1/n}.$$\eL
\begin{pf} Fix $t$ sufficiently large so that
$E_t=\{G\le-t\}\sbs D_{r/2}$. Then $t\om^*({\bf z},E_t,D)=G_t({\bf
z})=\max\{G({\bf z}),-t\}$ and from \cite[Theorem 4.2]{De}
\begin{equation}\label{d:si}
\il D(dd^cG)^n=\il D(dd^cG_t)^n.
\end{equation}
If $E$ is relatively compact in $D$, then by \cite[Prop.
4.7.2]{K}$$C^*(E)=\il D(dd^c\om^*({\bf z},E,D))^n.$$ Hence by
(\ref{e:si}) and (\ref{d:si})
\begin{equation}\begin{align}
&\il D|G_t|(dd^c\rho)^n=t\il D|\om^*({\bf
z},E_t,D)|(dd^c\rho)^n\le \notag\\&t\beta(r/2)\left(\il
D(dd^c\om^*({\bf z},E_t,D))^n\right)^{1/n}=\beta(r/2)\left(\il
D(dd^c G)^n\right)^{1/n}.\notag\end{align}\end{equation} By
monotone convergence
$$\il D|G|(dd^c\rho)^n=\lim_{t\to\infty}\il D|G_t|(dd^c\rho)^n,$$
and the lemma is proved.
\end{pf}
\par Now we can state a criterion for pluripolarity.
\bT\label{C:cpp} Under the hypotheses of Lemma \ref{L:ma} suppose
that $\{G_j\}$ is a sequence of multipole Green functions on $D$
with poles in $D_r$, $r<0$, and
$$\sup_j \il D(dd^cG_j)^n\le A<\infty.$$ For $t>0, \ s<0$ let
$$E_t=\{{\bf z}\in D_s:\limsup_{j\to \infty} G_j({\bf z})<-t\}.$$ Then for
the outer capacity of $E_t$ relative to $D$ we have the estimate
$$C^*(E_t)\le\frac{b(r)}{t\al(s)}A^{1/n}$$
where $\al(s)$ and $b(r)$ are positive constants depending only on
$s$ and $r$. In particular, the set
$$E=\{{\bf z}\in D:\lim_{j\to \infty}
G_j({\bf z})=-\infty\}$$ is pluripolar. \eT
\begin{pf} Define
$$E^k_t=\{{\bf z}\in E_t:G_j({\bf z})<-t \ \hbox{for all} \ j\ge k\}.$$
Let $\om^*$ be the relative extremal function of $E_t^k$ in $D$.
Then $G_j\le t\om^*$ for all $j\ge k$. Therefore, by (\ref{e:si})
$$\il D|G_j|(dd^c\rho)^n\ge t\il D|\om^*|(dd^c\rho)^n\ge
t\al(s)C^*(E_t^k).$$ By Lemma \ref{L:ma}
$$\il D|G_j|(dd^c\rho)^n\le b(r) A^{1/n}.$$
Thus $$C^*(E_t^k)\le \frac{b(r)}{t\al(s)}A^{1/n}.$$
\par Since $E^k_t\subset E^{k+1}_t$ and $E_t=\cup_kE^k_t$, by
\cite[Corollary 4.7.11]{K}
$$C^*(E_t)=\lim_{k\to\infty}C^*(E^k_t).$$
Hence we get the first statement of our lemma.
\par It follows from this statement that for every $s<0$ the set
$E\cap D_s$ is pluripolar. Hence $E$ is also pluripolar.
\end{pf}
\par We include the following lemma  because we were not able to find it
in the literature. We denote by $g_D(z,w)$ the (negative) Green
function of a domain $D\sbs{\Bbb C}$ with pole at $w\in D$.
\bL\label{L:gfr} Fix $r, a$ with  $r>a>1$ and let $A=\{z\in{\Bbb
C}:\,1/r<|z|<r\}$ be an annulus. Let $g(z)=g_A(z,1)$ be the Green
function of $A$ with pole at $w=1$. Then
$$\sup_{|z|=1}g(z)\le c(a)\log r,$$ where $c(a)<0$ depends only on
$a$.
\eL
\begin{pf} The mapping
$$w=f_1(z)=\exp\left(\log r\left(1+i\frac2\pi z\right)\right)$$
is a universal holomorphic covering map of $A$ by the strip
$T=\{z\in{\Bbb C}:\,0<\hbox{\bf Im} \ z<\pi\}$. Since $f_1(z)=1$
if and only if $$z=z_k=\frac{k\pi^2}{\log r}+\frac\pi2i, \
k=0,\pm1,\pm2,...$$ we see that
$$g_A(f_1(z),1)=\sum_{k=-\infty}^\infty g_T(z,z_k).$$
Note that if $I_1=\{x+i\pi/2:\,-\pi^2/\log r \le x\le0\}$, then
$f_1(I_1)=\{|z|=1\}$.
\par The mapping
$$f_2(z)=\log\left(i\frac{1-z}{1+z}\right)$$
maps the unit disk $U$ conformally onto the strip $T$. Since
$f_2(z'_k)=z_k$ and $f_2(I_2)=I_1$, where
$$z'_k=\frac{1-e^{k\pi^2/\log r}}{1+e^{k\pi^2/\log r}}$$
and $$I_2= \left[0,\frac{1-e^{-\pi^2/\log r}}{1+e^{-\pi^2/\log
r}}\right],$$ we see that to prove the lemma we need to estimate
the supremum $M$ of
$$u(x)=\sum_{k=-\infty}^\infty g_U(x,z'_k)=
\sum_{k=-\infty}^\infty\log\left|\frac{x-z_k'}{1-z_k'x}\right|$$
when $x$ varies over $I_2$.

If $k\le-1$, then the maximum of the function
$\log|(x-z'_k)/(1-z'_kx)|$ is achieved when $x=0$. Thus
$$M\le
\sum_{k=-\infty}^{-1} \log\frac{1-e^{k\pi^2/\log
r}}{1+e^{k\pi^2/\log r}}=\sum_{k=-\infty}^{-1}
\log\left(1-\frac{2e^{k\pi^2/\log r}}{1+e^{k\pi^2/\log
r}}\right).$$ Since $\log(1-x)\le -x$,
$$M\le-\sum_{k=-\infty}^{-1}
\frac{2e^{k\pi^2/\log r}}{1+e^{k\pi^2/\log r}}\le
-\sum_{k=-\infty}^{-1}e^{k\pi^2/\log r}=-\frac{e^{-\pi^2/\log
r}}{1-e^{-\pi^2/\log r}}.$$ Since $1-e^{-x}\le x$, we finally
obtain
$$M\le-\frac{e^{-\pi^2/\log r}}{\pi^2}\log r\le
-\frac{e^{-\pi^2/\log a}}{\pi^2}\log r.$$
\end{pf}
\par The corollary below allows us to verify the conditions of Theorem
\ref{C:cpp} by constructing special sequences of holomorphic
mappings of appropriate annuli.
\bC\label{P:annuli}
Let $E$ be a set in a ball $B\sbs{\Bbb C}^n$. Suppose that there
exists a sequence of arrays of points
$W_k=\{w_{1k},\dots,w_{m_kk}\}$ in $B$ and a sequence of numbers
$r_k$ such that, for each ${\bf z}\in E$, there exists a sequence
of holomorphic mappings $f_k$ of annuli $A_k=\{\zeta \in {\Bbb
C}:1/r_k<|\zeta|<r_k\}$ into $B$ satisfying:
\be
\item $f_k(e^{2\pi ij/m_k})=w_{jk}$, $1\le j\le m_k$;
\item there exists $\zeta_k\in S$ with $f_k(\zeta_k)={\bf z}$;
\item $$\limsup_{k\to\infty}m_k^{1-1/n}\log r_k=\infty.$$
\ee
Then the set $E$ is pluripolar.
\eC
\begin{pf} By passing to a subsequence we may assume that
$$\lim_{k\to\infty}m_k^{1-1/n}\log r_k=\infty.$$
Let $B'$ be the ball concentric with $B$ and of twice the radius. We
let $G_k({\bf z})$ denote the multipole Green function on $B'$
with poles at the points in $W_k$ of weight $m_k^{-1/n}$. Then
$$\il {B'}(dd^cG_k)^n\le A<\infty.$$
For ${\bf z}\in E$ we let $u_k(\zeta)=G_k(f_k(\zeta))$ so that
$u_k(\zeta_k)=G_k({\bf z})$. The functions $u_k(\zeta)$ are
negative and subharmonic and have poles of order $m_k^{-1/n}$ at
the points $e^{2\pi ij/m_k}$, $1\le j\le m_k$. Letting
$g(t)=g_R(t,1)$ denote the Green function of the annulus $R:=\{t
\in {\Bbb C}:r_k^{-m_k}<|t|<r_k^{m_k}\}$ with pole at $t =1$, we
have $u_k(\zeta)\le m_k^{-1/n}g(\zeta^{m_k})$. Note that
$r_k^{m_k}\to \infty$; thus for $k$ sufficiently large,
$r_k^{m_k}>2$. By Lemma \ref{L:gfr} with $a=2$, for such $k$
$$G_k({\bf z})=u_k(\zeta_k)\le c(2)m_k^{-1/n}\log
r_k^{m_k}=c(2)m_k^{1-1/n}\log r_k,$$ where $c(2)<0$. Thus
$$\lim_{k\to\infty}G_k({\bf z})=-\infty.$$ By Theorem
\ref{C:cpp} the set $E$ is pluripolar. \end{pf}

\section{Quasianalytic functions in the sense of Denjoy}\label{S:dq}
\par The proofs of Theorem \ref{T:bq} for Denjoy quasianalytic
functions  and of Theorem \ref{T:dq} will follow from a slightly
more general, albeit slightly technical-looking, result. Given a
smooth ($C^\infty$) function $f$ on $S$ with associated function
$\tau_f(r)$, we define
\begin{eqnarray}\label{e:tn}
\log  t_n=\min\left\{-\frac{\log r^3\tau_f(r)}{r}:\, 1\leq r\leq
n\right\}.\end{eqnarray}
\begin{Proposition}\label{P:dq} Let $f:S\ra{\Bbb C}$ be a
$C^\infty$ function such that
\begin{eqnarray}\label{e:hyp}
 & & \limsup_{n\ra\infty}\sqrt{n}\,\log  t_n=\infty.
\end{eqnarray}
Then the set $\Gamma_f$ is pluripolar.
\end{Proposition}
\begin{pf} Let $f:S\ra{\Bbb C}$ be a smooth function with
Fourier expansion (\ref{e:Fourier}). By (\ref{e:deriv}) we have
\begin{equation}\label{e:coef}
|c_k|\leq M_j(f)/|k|^j,\;|k|\geq1,\;j\geq0.
\end{equation}
\par The idea of the proof is to interpolate $f$ by trigonometric
polynomials $L_n(f,z_0;z)$ at the $n$-th roots of unity and at
some other point $z_0$ with $|z_0|=1$. Lemma \ref{L:er} provides
estimates for $L_n(f,z_0;z)$ when $z$ lies in an annulus
$A(t)=\{z:\,1/t\leq|z|\leq t\}$. It follows from these estimates
that we may apply Corollary \ref{P:annuli} using the sequence of
arrays of points $W_n=\{(z,f(z)):\,z^n=1\}$, the annuli $A(t_n)$
and our hypothesis (\ref{e:hyp}).
\par We proceed to define the appropriate interpolating trigonometric
polynomials $\{L_n(f,z_0;z)\}$; recall
these are rational functions restricted to $S$. To begin with, we let
$$L_n(f;z)=\sum_{r=0}^{n-1}a_{n,r}z^r+
\sum_{r=1}^{n}\frac{b_{n,r}}{z^r},$$ where
$$a_{n,r}=\sum_{j=0}^\infty c_{r+nj},\,
     b_{n,r}=\sum_{j=0}^\infty c_{-r-nj},$$ and then we define
$$L_n(f,z_0;z)=L_n(f;z)+\frac{z^n-1}{z_0^n-1}\,
(f(z_0)-L_n(f;z_0)).$$ Here $z_0$ is any point on the unit circle
with $z_0^n\neq1$. Note that $a_{n,r},\,b_{n,r}$ are well
defined since $\sum|c_k|<\infty$ from (\ref{e:coef}).
\par If $z_1,\dots,z_n$ are the $n$-th roots of unity, then
$z_l^{r+nj}=z_l^r$ and, therefore,
$$L_n(f,z_0;z_l)=f(z_l),\;l=0,1,\dots,n,$$ i.e., $L_n(f,z_0;z)$
interpolates $f$ at $z_0,\dots,z_n$. If $z_0^n=1$ we define
$L_n(f,z_0;z)=L_n(f;z)$. We remark that $C_f$ will denote a
constant depending on $f$ which may vary from line to line.

\begin{Lemma}\label{L:er}
There exists a constant $C_f$ depending only on $f$ so that for
every $n\geq1$, every $z_0$ with $|z_0|=1$, and every $t>1$, we
have
$$|L_n(f,z_0;z)|\leq C_f\left(1+\sum_{r=1}^n r\tau_f(r)t^r\right)$$
for all $z$ with $1/t\leq|z|\leq t$.
\end{Lemma}

\begin{pf} If $s\geq2$ we have, using (\ref{e:coef}),
$$|a_{n,r}|\leq\sum_{j=0}^\infty|c_{r+nj}|\leq
M_s(f)\sum_{j=0}^\infty\frac{1}{(r+nj)^s}\leq
2\,\frac{M_s(f)}{r^s}\;,$$ and similarly
$|b_{n,r}|\leq2M_s(f)/r^s$. By the definition of $\tau_f(r)$ we
clearly have
$$\tau_f(r)=\min\left\{M_0(f),\frac{M_1(f)}{r},\frac{M_2(f)}{r^2},
\inf_{s\geq3}\frac{M_s(f)}{r^s}\right\}
=\inf_{s\geq3}\frac{M_s(f)}{r^s}$$ for all $r>r_0(f)$. We conclude
that $|a_{n,r}|,\,|b_{n,r}|$ are bounded above by $C_f\tau_f(r)$.
Therefore
$$|L_n(f;z)|\leq C_f\left(1+2\sum_{r=1}^n\tau_f(r)t^r\right)$$
holds for $1/t\leq|z|\leq t$. This gives the desired bound when
$z_0^n=1$.
\par Next, if $z_0^n\neq1$, we note that
$$|(z_0^{nj}-1)/(z_0^n-1)|=|z_0^{n(j-1)}+z_0^{n(j-2)}+\dots+1|\leq j,$$
moreover,
$$f(z_0)=\sum_{j=-\infty}^{\infty}c_jz_0^j=
\sum_{r=0}^{n-1}\sum_{j=0}^\infty c_{r+nj}z_0^{r+nj}
+\sum_{r=1}^n\sum_{j=0}^\infty \frac{c_{-r-nj}}{z_0^{r+nj}},$$
thus
\begin{eqnarray*}
 & & \left|\frac{f(z_0)-L_n(f;z_0)}{z_0^n-1}\right|=\\
 & & \left|\sum_{r=0}^{n-1}\sum_{j=0}^\infty c_{r+nj}\,
\frac{z_0^{r+nj}-z_0^r}{z_0^n-1}+ \sum_{r=1}^n\sum_{j=0}^\infty
\frac{c_{-r-nj}(1-z_0^{nj})}{z_0^{r+nj}(z_0^n-1)}\right|\leq\\
 & & \sum_{r=0}^{n-1}\sum_{j=0}^\infty j|c_{r+nj}|+
\sum_{r=1}^n\sum_{j=0}^\infty j|c_{-r-nj}|.
\end{eqnarray*}
For $s\geq3$ we have, using (\ref{e:coef}),
$$\sum_{j=1}^\infty
j|c_{\pm(r+nj)}|\leq\sum_{j=1}^\infty\frac{jM_s(f)}{(r+nj)^s}
\leq\sum_{j=1}^\infty\frac{jM_s(f)}{n^sj^s}
\leq2\,\frac{M_s(f)}{n^s}\;.$$ Thus by the definition of $\tau_f(n)$,
$$\left|\frac{f(z_0)-L_n(f;z_0)}{z_0^n-1}\right|\leq4n\tau_f(n).$$
Since $|z^n-1|\leq2t^n$ for $|z|\leq t$, this estimate, together
with the bound on $|L_n(f;z)|$, implies the lemma. \end{pf}

\par Using (\ref{e:tn}) we have defined the sequence $\{t_n\}$ via
$$t_n =\min \left\{ \frac{1}{(r^3\tau_f(r))^{1/r}}:1\leq r \leq
n\right\}.$$
Since the numbers $t_n$ are decreasing, we obtain from
(\ref{e:hyp}) that $t_n>1$ for each $n$; moreover
$r^3\tau_f(r)t_n^r\leq1$ for $r\leq n$. If $1/t_n\leq|z|\leq t_n$
it follows from the previous lemma that
$$|L_n(f,z_0;z)|\leq
C_f\left(1+\sum_{r=1}^n\frac{r^3\tau_f(r)t_n^r}{r^2}\right)
\leq C_f\left(1+\sum_{r=1}^n\frac{1}{r^2}\right)\leq 3C_f.$$

\par We conclude that for each $n$ and for each $z_0$ with $|z_0|=1$,
the images of the annuli $A(t_n)=\{1/t_n\le|z|\le t_n\}$ under
the mappings $h_n(z)=(z,L_n(f,z_0;z))$ are contained in a ball
$B\subset{\Bbb C}^2$ centered at the origin and of radius $R_f$
depending only on $f$. Thus by (\ref{e:hyp}) and Corollary
\ref{P:annuli} the set $\Gm_f$ is pluripolar.
\end{pf}

\par Before proceeding with the proofs of Theorems \ref{T:bq}
(for functions quasianalytic in the sense of Denjoy) and
\ref{T:dq} we make the following remarks. Note that the graph of
$f$ is pluripolar if and only if the graph of $cf$ is pluripolar,
where $c\neq0$ is a constant. Multiplying $f$ by a small constant
we may assume in the sequel that our smooth functions $f$ verify
$M_3(f)<1/2$. Let
\begin{eqnarray}\label{e:theta}
\wtl \tau_f(r):=\inf_{s\geq 3}\frac{M_s(f)}{r^{s-3}}=\inf_{s\geq
0}\frac{M_{s+3}(f)}{r^{s}}<\frac{1}{2}
\end{eqnarray}
be the associated function for the shifted sequence $\{\tilde
M_s\}=\{M_{s+3}(f)\}$; setting
\begin{eqnarray}\label{e:then}
\log \th_f(n)=\min\left\{-\frac{\log\wtl\tau_f(r)}r:\,1\leq r\leq
n\right\},
\end{eqnarray}
it follows that
$$\log t_n \geq \log \th_f(n)>0.$$
\par
From the definition of $\tau_f(r)$, as indicated in the proof of
Lemma \ref{L:er}, we clearly have
$$r^3\tau_f(r)=
\min\left\{M_0(f)r^3,M_1(f)r^2,M_2(f)r,
\inf_{s\geq3}\frac{M_s(f)}{r^{s-3}}\right\}=\wtl\tau_f(r)$$ for
all $r>r_0(f)\ge0$. We show that if $f$ is quasianalytic in the
sense of Denjoy or if $f$ belongs to a Gevrey class $G\{L_j\}$,
where $L_j=o(j^2)$, then
\begin{equation}\label{e:thn}
\limsup_{n\to\infty}\sqrt n \log \th_f(n)= \infty.
\end{equation}
This implies condition (\ref{e:hyp}) from Proposition \ref{P:dq}
and finishes the proofs of Theorems \ref{T:bq} and \ref{T:dq}.

\par We need the following lemma, whose proof we postpone until the end
of this section.

\begin{Lemma}\label{L:calc}
Let $\tilde h(s)=h(e^s)$ be a positive, increasing, convex
function of $s$ on $[0,\infty)$ and let $\tilde H(x)=\min\{\tilde
h(s)e^{-s}:\,0\le s\le x\}$. If $\tilde H(x)\le Ce^{-x/2}$ for all
$x\ge 0$, then
$$\il 0^\infty \tilde h(s)e^{-s}\,ds =\il 1^\infty
\frac{h(t)}{t^2}dt<\infty.$$
\end{Lemma}

\par We now prove Theorem \ref{T:bq} for quasianalytic functions in the
sense of Denjoy.
\begin{pf} The condition (\ref{Q:def}) for quasianalyticity holds
for $\wtl\tau_f(r)$ (equivalently, for $\{\tilde M_s\}$); in
particular, since $r^3\tau_f(r)=\wtl\tau_f(r)$ for all $r>r_0(f)$,
using (\ref{Q:def}) for $\tau_f(r)$, we have
$$\int_{r_0}^{\infty}\frac{\log
\wtl\tau_f(r)}{1+r^2}\,dr=\int_{r_0}^{\infty}\frac{3\log
r}{1+r^2}\,dr+ \int_{r_0}^{\infty}\frac{\log
\tau_f(r)}{1+r^2}\,dr=-\infty;$$ i.e., letting $h(r):=-\log
\wtl\tau_f(r)$,
\begin{equation}\label{ue:unint}
\il 1^\infty \frac{h(t)}{1+t^2}dt= \infty.
\end{equation}
We assume, for the sake of obtaining a contradiction, that (\ref{e:thn})
does not hold. Then, with
$h(t)=-\log \wtl\tau_f(t)$ and $\tilde h(s)=h(e^s)$,
we have
$$H(t)=\log \th_f(t)=\min\{h(r)/r:\,1\leq r\leq t\},$$
and it follows that there exists a constant $C>0$ so that
$H(n)<C/\sqrt{n}$ for every integer $n$. Since $H$ is decreasing,
$H(x)\le2C/\sqrt{x}$ for all $x$. In terms of
$$\tilde H(x)=\min\{\tilde h(s)e^{-s}:\,0\le s\le x\}=H(e^x),$$
$$\tilde H(x)<2Ce^{-x/2}\,,\;\forall\,x\geq0.$$
By (\ref{e:theta}) we have $h(t)>0$, so Lemma \ref{L:calc} implies
that
\begin{equation}\label{e:integral}
\int_1^\infty \frac{h(t)}{t^2}<\infty,
\end{equation}
which contradicts (\ref{ue:unint}); i.e., the fact that $\wtl
\tau_f(r)$ (equivalently, $\{\tilde M_s\}$) defines a
quasianalytic class.
\end{pf}
\par We proceed with the proof of Theorem \ref{T:dq}.
\begin{pf} Writing $L_j=j^2/a_j$, $a_j\ra\infty$, we see that $f\in
G\{L_j\}$ satisfies
$$M_j(f)/r^j\leq(Cj^2/(ra_j))^j$$
for all $j,r\geq1$. Taking $j=[{\sqrt{r}}]$, the greatest integer
in $\sqrt{r}$, we obtain
$$\wtl\tau_f(r)\leq r^3(C/a_{[\sqrt{r}]})^{[\sqrt{r}]};$$
hence
$$\frac{-\log\wtl\tau_f(r)}r\geq
-3\frac{\log r}r+\frac{[\sqrt{r}]}r\, \left(\log
a_{[\sqrt{r}]}-\log C\right).$$ Since the nonnegative function $\log
\th_f(r)$
is decreasing, either $\log \th_f(r)=c>0$ for
large $r$ or else there is a sequence $r_k\to\infty$ such that
$$\log \th_f(r_k)=-\frac{\log\wtl\tau_f(r_k)}{r_k}.$$
In the first case the condition (\ref{e:thn}) is clearly
satisfied. In the second case, if $n_k=[r_k]-1$, then
$$\sqrt{n_k}\log \th_f(n_k)\ge-3\frac{\sqrt{n_k}\log r_k}{r_k}+
\frac{\sqrt{n_k}[\sqrt{r_k}]}{r_k}\, \left(\log
a_{[\sqrt{r_k}]}-\log C\right)\to\infty.$$
\end{pf}

{\bf Proof of Lemma \ref{L:calc}}. Since the functions $\tilde h$
and $\tilde H$ are continuous, the set $E=\{x:\,\tilde
h(x)e^{-x}=\tilde H(x)\}$ is closed. We have
$$\il E \tilde h(x)e^{-x}\,dx=\int_E\tilde H(x)\,dx\leq
C\int_0^\infty e^{-x/2}<\infty.$$
Let $F=[0,\infty]\sm E$. The set $F$ is open and, therefore,
$F=\cup(a_j,b_j)$, where $a_j,b_j\in E$ or $b_j=\infty$.
\par We first show that $\tilde H(x)=\tilde h(a_j)e^{-a_j}$ on
$[a_j,b_j]$. Indeed,
since $\tilde H(x)$ is decreasing, $\tilde H(x)\le \tilde h(a_j)e^{-a_j}$
on
$[a_j,b_j]$. If $\tilde h(x)e^{-x}< \tilde h(a_j)e^{-a_j}$ for some
$x\in[a_j,b_j]$, then let $s$ be a point in $[a_j,x]$, where
$\tilde h(x)e^{-x}$ attains its minimum. Clearly $s>a_j$. But then
$\tilde H(s)=\tilde h(s)e^{-s}$ and $s\in E$. Thus $s=b_j$. Hence $\tilde
h(x)e^{-x}\ge
\tilde h(a_j)e^{-a_j}$ for all $x\in [a_j,b_j]$ and
$\tilde H(x)=\tilde h(a_j)e^{-a_j}$.
\par Next we show that $b_j<\infty$ for all $j$. If not, then
$\tilde H(x)=\tilde h(a_j)e^{-a_j}\le Ce^{-x/2}$ for all $x\ge
a_j$ and hence $\tilde h(a_j)=0$. This contradiction shows that
$b_j<\infty$.
\par From the above results we have
$$\tilde H(x)=\tilde h(a_j)e^{-a_j}=\tilde h(b_j)e^{-b_j}\le Ce^{-b_j/2}$$
when
$a_j\le x\le b_j$ Consequently, $\tilde h(a_j)\le Ce^{a_j-b_j/2}$ and
$\tilde h(b_j)\le Ce^{b_j/2}$.
\par If $x=\al a_j+(1-\al) b_j$, $0\le\al\le1$, then
$$\tilde h(x)\le\al \tilde h(a_j)+(1-\al)\tilde h(b_j)\le
Ce^{-b_j/2}\left(\al e^{a_j}+(1-\al)e^{b_j}\right)$$ and
$$\tilde h(x)e^{-x}\le Ce^{-b_j/2}\frac{\al e^{a_j}+(1-\al)e^{b_j}}
{e^{\al a_j+(1-\al) b_j}}=Ce^{-b_j/2}\frac{\al
e^{c_j}+1-\al}{e^{\al c_j}},$$where $c_j=a_j-b_j$.
\par We split the intervals $(a_j,b_j)$ into three separate types.
The first type consists of all intervals having length at most
$2$; i.e., $-c_j\le2$. On these intervals
$$\tilde h(x)e^{-x}\le Ce^2e^{-b_j/2}\le Ce^2e^{-x/2}$$ for
$x\in [a_j,b_j]$. Letting $F_1$ denote the union of these type one
intervals, we have
$$\il {F_1} \tilde h(x)e^{-x}\,dx<\infty.$$
\par The second type of interval is one of the form $[a_j,b_j]$ where
$b_j-a_j>2$ and for which $\tilde h(x)\le e^{x/2}$ for $a_j\le
x\le b_j$. If $F_2$ is the union of these intervals, then
$$\il {F_2}\tilde h(x)e^{-x}\,dx\le\il 0^{\infty}e^{-x/2}\,dx<\infty.$$
\par We are left with type three intervals $[a_j,b_j]$ where $-c_j\ge2$
and there exists a point $x_j$ between $a_j$ and $b_j$ with
$\tilde h(x_j)\ge e^{x_j/2}$. If $x=\al a_j+(1-\al) b_j$,
$0\le\al\le1$, then $x=c_j\al+b_j$. Hence
\begin{equation}\begin{align}
&\il {a_j}^{b_j}\tilde h(x)e^{-x}\,dx\le -Ce^{-b_j/2}c_j\il 0^1\left(\al
e^{-\al
c_j}e^{c_j}+(1-\al)e^{-\al c_j}\right)\,d\al=\notag\\
&Ce^{-b_j/2}\frac{2-(e^{c_j}+e^{-c_j})}{c_j}\le
-Ce^{-b_j/2}\frac{e^{c_j}+e^{-c_j}}{c_j}.\notag
\end{align}\end{equation}
We enumerate these intervals consecutively. Recall that $c_j<-2$;
hence $b_{j+1}>b_j+2$, so $b_j>2j$. Thus
$$-\sum_je^{-b_j/2}\frac{e^{c_j}}{c_j}<\infty.$$ Since
$\tilde h(a_j)e^{-a_j}\le Ce^{-b_j/2}$, $e^{b_j/2-a_j}\le C/\tilde
h(a_j)$.
Hence
$$-Ce^{-b_j/2}\frac{e^{-c_j}}{c_j}=C\frac{e^{b_j/2-a_j}}{b_j-a_j}\le
C^2\frac1{\tilde h(a_j)(b_j-a_j)}.$$
\par But $\tilde h(a_{j+1})\ge \tilde h(x_j)\ge e^{a_j/2}$. Since
$a_{j+1}\ge
2j$, we see that
$$-\sum_je^{-b_j/2}\frac{e^{-c_j}}{c_j}<\infty;$$
hence
$$\il 0^\infty \tilde h(x)e^{-x}\,dx<\infty.$$
$\Box$
\par {\bf Final remark.} For a $C^\infty$
mapping $f=(f_1,\dots,f_N):\,S\to {\Bbb C}^N$ we can define
$M_j(f)=\sup \{M_j(f_k):\,1\le k\le N\}$. Referring to Section
\ref{S:bdf}, we can define the associated function $\tau_f(r)$ and
the sequence $\{t_n\}$ as well as Gevrey classes $G\{L_j\}$ of
mappings. One can verify that the graph $\Gm_f$ of $f$ in ${\Bbb
C}^{N+1}$ is pluripolar if condition (\ref{e:hyp}) is replaced by
$$\limsup_{n\to\infty}n^{1-1/(N+1)}\log t_n=\infty.$$
Moreover, an argument similar to one used in the proof of Theorem
\ref{T:dq}
can be used to prove that if $f\in G\{L_j\}$, where
$L_j=o(j^{N+1})$, then $\Gm_f$ is pluripolar.

\end{document}